\documentclass[11pt]{amsart}

\usepackage{amsmath,amssymb,amsthm,mathtools}
\usepackage[margin=1.15in]{geometry}
\usepackage{enumitem}
\usepackage{booktabs}
\usepackage{microtype}
\usepackage[hidelinks]{hyperref}

\theoremstyle{plain}
\newtheorem{theorem}{Theorem}[section]
\newtheorem{proposition}[theorem]{Proposition}
\newtheorem{lemma}[theorem]{Lemma}
\newtheorem{corollary}[theorem]{Corollary}
\theoremstyle{definition}

\newtheorem{example}[theorem]{Example}
\newtheorem{remark}[theorem]{Remark}
\newtheorem{question}[theorem]{Question}

\numberwithin{equation}{section}
\newcommand{\R}{\mathbb{R}}
\newcommand{\Z}{\mathbb{Z}}
\newcommand{\N}{\mathbb{N}}
\newcommand{\T}{\mathbb{T}}
\newcommand{\Et}{\widetilde{E}}
\newcommand{\Ct}{\widetilde{C}}
\newcommand{\Mt}{\widetilde{M}}
\newcommand{\cB}{\mathcal{B}}
\newcommand{\BHL}{\mathcal{B}_{\mathrm{HL}}}
\newcommand{\BS}{\mathcal{B}_{\mathrm{S}}}
\newcommand{\CHL}{C_{\mathrm{HL}}}

\newcommand{\CtHL}{\Ct_{\mathrm{HL}}}
\newcommand{\CsHL}{C^{*}_{\mathrm{HL}}}
\newcommand{\MsHL}{M^{*}_{\mathrm{HL}}}
\newcommand{\MtHL}{\Mt_{\mathrm{HL}}}
\newcommand{\CS}{C_{\mathrm{S}}}
\newcommand{\CtS}{\Ct_{\mathrm{S}}}
\newcommand{\CsS}{C^{*}_{\mathrm{S}}}
\newcommand{\Lat}[1]{(#1)_{\Z}}
\newcommand{\one}{\mathbf{1}}

\begin{document}

\title[Sharp Solyanik estimates for ergodic maximal operators]{Sharp Solyanik estimates for discrete and ergodic\\ Hardy--Littlewood maximal operators}

\author{Kesen Huang}
\thanks{This paper was prepared with substantial assistance from Claude (Anthropic), which proposed the main arguments. The author completed the gaps in the proofs and independently verified all results; any remaining errors are the author's.}
\address{Department of Mathematics, Baylor University, Waco, Texas}
\email{Kesen\_Huang1@baylor.edu}

\subjclass[2020]{Primary 37A46, 42B25; Secondary 37A30, 11P21}
\keywords{Solyanik estimates, Tauberian constants, halo functions, maximal operators, ergodic theory, Rokhlin towers, lattice points}

\begin{abstract}
Let $U_1,\dots,U_n$ be commuting invertible measure-preserving transformations of a probability space, and let $\CsHL(\alpha)$ denote the sharp Tauberian constant of the associated ergodic uncentered Hardy--Littlewood maximal operator. Hagelstein and Parissis showed that $\CsHL(\alpha)-1\lesssim_n(1/\alpha-1)^{1/(n(n+1))}$ as $\alpha\to1^-$ and asked whether this exponent is sharp. We prove that $\CsHL(\alpha)-1\lesssim_n(1/\alpha-1)^{1/n}$ for all such transformations, and that this is sharp: equality of the ergodic and discrete Tauberian constants holds for every action admitting arbitrarily large Rokhlin towers (in particular for every free action, ergodic or not), and the discrete constant satisfies $\CtHL(\alpha)-1\asymp_n(1/\alpha-1)^{1/n}$. For $n\ge2$ the sharp ergodic exponent $1/n$ is strictly smaller than the sharp Euclidean exponent $2/(n+1)$ recently obtained by Mukherjee; the discrepancy is caused by a lattice ``first layer'' effect. We also show that the discrete Tauberian constants of the Hardy--Littlewood and strong maximal operators dominate their Euclidean counterparts; for the strong maximal operator this yields $\CsS=\CS$ for every free action.
\end{abstract}

\maketitle

\section{Introduction}

\subsection{Solyanik estimates}

By a \emph{basis} we mean a collection $\cB$ of bounded open subsets of $\R^n$ containing the origin. Following \cite{HP16}, we associate with $\cB$ three maximal operators. The geometric maximal operator acts on $f\in L^1_{\mathrm{loc}}(\R^n)$ by
\[
M_\cB f(x)\coloneqq\sup_{R\in\cB}\frac{1}{|R|}\int_R|f(x+y)|\,dy ,
\]
the discrete maximal operator acts on $f\colon\Z^n\to\R$ by
\[
\Mt_\cB f(m)\coloneqq\sup_{R\in\cB}\frac{1}{\#(R\cap\Z^n)}\sum_{j\in R\cap\Z^n}|f(m+j)| ,
\]
and, given commuting invertible measure-preserving transformations $U_1,\dots,U_n$ of a probability space $(\Omega,\Sigma,\mu)$, the ergodic maximal operator acts on $f\in L^1(\Omega)$ by
\[
M^*_\cB f(\omega)\coloneqq\sup_{R\in\cB}\frac{1}{\#(R\cap\Z^n)}\sum_{j\in R\cap\Z^n}|f(U^j\omega)| ,\qquad U^j\coloneqq U_1^{j_1}\cdots U_n^{j_n} .
\]
For $0<\alpha<1$ the corresponding \emph{sharp Tauberian constants} are
\begin{gather*}
C_\cB(\alpha)\coloneqq\sup_{0<|E|<\infty}\frac{|\{M_\cB\one_E>\alpha\}|}{|E|},\qquad
\Ct_\cB(\alpha)\coloneqq\sup_{0<\#\Et<\infty}\frac{\#\{\Mt_\cB\one_{\Et}>\alpha\}}{\#\Et},\\
C^*_\cB(\alpha)\coloneqq\sup_{\mu(E)>0}\frac{\mu(\{M^*_\cB\one_E>\alpha\})}{\mu(E)} .
\end{gather*}
The function $\alpha\mapsto C_\cB(\alpha)$ is also called the \emph{halo function} of $\cB$. We are concerned with the bases $\BHL$ of all open Euclidean balls containing the origin and $\BS$ of all open axis-parallel rectangles containing the origin, and write $\CHL,\CtHL,\CsHL$ and $\CS,\CtS,\CsS$ for the associated constants. Since some elements of these bases contain lattice points with negative coordinates, the definition of $M^*_\cB$ requires the transformations to be invertible; this assumption is implicit in \cite{HP16}. No ergodicity assumption is made anywhere in this paper.

Solyanik~\cite{So} proved that $\CHL(\alpha)\to1$ and $\CS(\alpha)\to1$ as $\alpha\to1^-$. Quantitative versions of such statements are called \emph{Solyanik estimates}. Hagelstein and Parissis~\cite{HP14} showed that
\[
\CHL(\alpha)-1\lesssim_n(1/\alpha-1)^{1/(n+1)},\qquad \CS(\alpha)-1\lesssim_n(1/\alpha-1)^{1/n}
\]
for $\alpha$ close to $1$, and used Solyanik estimates in \cite{HP15} to prove that $\CHL$ and $\CS$ are locally H\"older continuous on $(0,1)$. In \cite{HP16} they initiated the study of Solyanik estimates in ergodic theory. Combining Calder\'on's transference principle~\cite{Ca}, which gives $C^*_\cB\le\Ct_\cB$, with comparisons between discrete and Euclidean constants, they obtained
\begin{equation}\label{eq:HP16}
\CsS(\alpha)-1\lesssim_n(1/\alpha-1)^{1/n},\qquad \CsHL(\alpha)-1\lesssim_n(1/\alpha-1)^{1/(n(n+1))} .
\end{equation}
The second exponent arises as the product of two losses: the passage from lattice averages over balls to Euclidean averages costs a factor $1/n$ in the exponent, and the Euclidean estimate contributes $1/(n+1)$. Problem~7.2 of \cite{HP16} asks whether the exponent $1/(n(n+1))$ is sharp, and notes that it was likewise unknown whether $1/(n+1)$ is sharp for $\CHL$.

The Euclidean half of this problem was recently settled by Mukherjee~\cite{Mu}, who proved a sharp covering theorem for Euclidean balls and deduced that
\begin{equation}\label{eq:Mu}
\CHL(\alpha)-1\asymp_n(1-\alpha)^{2/(n+1)},\qquad n\ge2,
\end{equation}
confirming a conjecture of \cite{HP14}; see also \cite[Problem~8]{H24}. The purpose of the present paper is to settle the ergodic half. As we shall see, the answer is different from both $1/(n(n+1))$ and $2/(n+1)$.

Throughout, $A\lesssim_n B$ means $A\le C_nB$ for a constant $C_n$ depending only on $n$, and $A\asymp_nB$ means $A\lesssim_nB\lesssim_nA$. Since $1-\alpha\le1/\alpha-1\le2(1-\alpha)$ for $\alpha\ge1/2$, the quantities $1-\alpha$ and $1/\alpha-1$ may be used interchangeably in all asymptotic statements.

\subsection{Main results}

Our first result determines the Solyanik exponent of the discrete constant.

\begin{theorem}\label{thm:discrete}
For every $n\ge1$ there exist $c_n,C_n>0$ and $\alpha_n\in(0,1)$ such that
\[
c_n(1-\alpha)^{1/n}\le\CtHL(\alpha)-1\le C_n(1-\alpha)^{1/n},\qquad\alpha_n<\alpha<1 .
\]
\end{theorem}

The upper bound is proved directly on $\Z^n$, without passing through the Euclidean constant, using Mukherjee's covering theorem (for $n\ge3$ the earlier covering lemma of Han and Lu~\cite{HL} suffices; see Remark~\ref{rem:covering}). The lower bound comes from a simple example: a lattice cube of side length comparable to $(1-\alpha)^{-1/n}$, together with balls that are almost tangent to a lattice hyperplane and therefore pick up a single lattice point outside the cube.

Our second result is a two-sided transference principle for Tauberian constants. We say that the action $j\mapsto U^j$ \emph{admits arbitrarily large Rokhlin towers} if for every $K\in\N$ there is $\Omega_0\in\Sigma$ with $\mu(\Omega_0)>0$ such that the sets $U^j\Omega_0$, $j\in\{0,\dots,K-1\}^n$, are pairwise disjoint. Recall that the action is \emph{free} (or \emph{non-periodic} \cite{KW}) if $\mu(\{\omega:U^j\omega=\omega\})=0$ for every $j\ne0$.

\begin{theorem}\label{thm:transference}
Let $\cB$ be a basis in $\R^n$ and let $U_1,\dots,U_n$ be commuting invertible measure-preserving transformations of a probability space.
\begin{enumerate}[label=\textup{(\roman*)},leftmargin=2.2em]
\item \textup{(\cite[Theorem~3.1]{HP16})} $C^*_\cB(\alpha)\le\Ct_\cB(\alpha)$ for every $0<\alpha<1$.
\item If the action admits arbitrarily large Rokhlin towers, then $C^*_\cB(\alpha)=\Ct_\cB(\alpha)$ for every $0<\alpha<1$.
\end{enumerate}
On a standard probability space, the hypothesis in \textup{(ii)} holds whenever the set of points with trivial stabilizer has positive measure; in particular it holds for every free action, ergodic or not.
\end{theorem}

Combining the two theorems answers the ergodic part of \cite[Problem~7.2]{HP16}.

\begin{corollary}\label{cor:ergodic}
Let $n\ge1$, and let $c_n,C_n,\alpha_n$ be as in Theorem~\ref{thm:discrete}. For all commuting invertible measure-preserving transformations $U_1,\dots,U_n$ of a probability space and all $\alpha_n<\alpha<1$,
\[
\CsHL(\alpha)-1\le C_n(1-\alpha)^{1/n},
\]
and if the action admits arbitrarily large Rokhlin towers (for instance, if it is free), then also $\CsHL(\alpha)-1\ge c_n(1-\alpha)^{1/n}$. In particular, $1/n$ is the sharp Solyanik exponent for $\CsHL$ uniformly over all such transformations, and the exponent $1/(n(n+1))$ in \eqref{eq:HP16} is not sharp for any $n\ge1$.
\end{corollary}

Ergodicity plays no role here. In Example~\ref{ex:examples} we exhibit a free non-ergodic action for which $\CsHL=\CtHL$, and an ergodic action with bounded periods for which $\CsHL(\alpha)=1$ for all $\alpha$ sufficiently close to $1$. What governs the worst case is the availability of large Rokhlin towers.

For $n\ge2$ we have $1/n<2/(n+1)$, so by \eqref{eq:Mu} the discrete and ergodic constants converge to $1$ strictly more slowly than the Euclidean one. Our final result shows that the discrete constants always dominate the Euclidean ones for balls and rectangles.

\begin{theorem}\label{thm:comparison}
Let $\cB\in\{\BHL,\BS\}$. Then $C_\cB(\alpha)\le\Ct_\cB(\alpha)$ for every $0<\alpha<1$.
\end{theorem}

For the strong basis, the reverse inequality $\CtS\le\CS$ was proved in \cite[Section~4]{HP16}, using the fact that lattice averages over axis-parallel rectangles are Euclidean averages over axis-parallel rectangles. Combining it with Theorems~\ref{thm:comparison} and~\ref{thm:transference}, we obtain the following.

\begin{corollary}\label{cor:strong}
For every $0<\alpha<1$ we have $\CtS(\alpha)=\CS(\alpha)$. Consequently, $\CsS(\alpha)=\CS(\alpha)$ for every $0<\alpha<1$ and every action admitting arbitrarily large Rokhlin towers, in particular for every free action on a standard probability space.
\end{corollary}

Thus, for free actions the ergodic constant $\CsS$ inherits the local H\"older continuity of $\CS$ established in \cite{HP15}, which answers \cite[Problem~7.4]{HP16} for the strong maximal operator and free actions. The situation for balls is different: by Theorem~\ref{thm:discrete} and \eqref{eq:Mu}, for $n\ge2$ we have $\CHL(\alpha)<\CtHL(\alpha)$ for all $\alpha$ sufficiently close to $1$; Example~\ref{ex:strict} makes this explicit for $n=2$, where the strict inequality holds whenever $1-\alpha\le10^{-7}$. The following table summarizes the situation (for $n\ge2$; the ergodic column refers to actions admitting arbitrarily large Rokhlin towers).
\medskip
\begin{center}
\begin{tabular}{lll}
\toprule
& Euclidean constant & discrete constant $=$ ergodic constant (worst case)\\
\midrule
balls & exponent $2/(n+1)$ \ \cite{Mu} & exponent $1/n$ \ (Theorem~\ref{thm:discrete}, Corollary~\ref{cor:ergodic})\\
rectangles & $\CS$ & $\CtS=\CsS=\CS$ \ (Corollary~\ref{cor:strong})\\
\bottomrule
\end{tabular}
\end{center}
\medskip
The contrast between the two rows reflects curvature: a ball can touch a lattice hyperplane at a single lattice point, whereas an axis-parallel rectangle meets a layer of lattice points in a full $(n-1)$-dimensional block. We return to this point in Section~\ref{sec:remarks}.

\subsection{Organization}

Section~\ref{sec:prelim} collects notation and auxiliary results. The upper and lower bounds in Theorem~\ref{thm:discrete} are proved in Sections~\ref{sec:upper} and~\ref{sec:lower}. Section~\ref{sec:transference} contains the proofs of Theorem~\ref{thm:transference} and Corollary~\ref{cor:ergodic}, together with examples. Theorem~\ref{thm:comparison} and Corollary~\ref{cor:strong} are proved in Section~\ref{sec:comparison}, where we also give an explicit example of strict inequality for balls. Section~\ref{sec:remarks} contains further remarks and open questions.

\section{Preliminaries}\label{sec:prelim}

All balls are open. We write $B(c,r)$ for the open ball with center $c$ and radius $r$, with the convention $B(c,r)=\emptyset$ for $r\le0$, and $aB(c,r)\coloneqq B(c,ar)$ for $a>0$. For $A\subset\R^n$ let $\Lat{A}\coloneqq A\cap\Z^n$. Let $\kappa_n\coloneqq|B(0,1)|$ and $s_n\coloneqq\sqrt n/2$. For a finite set $\Et\subset\Z^n$ and $0<\alpha<1$ we write $\Et_\alpha\coloneqq\{m\in\Z^n:\MtHL\one_{\Et}(m)>\alpha\}$.

\begin{lemma}\label{lem:count}
For every finite family of balls $\{B(c_j,r_j)\}_j$ in $\R^n$,
\[
\Big|\bigcup_j B(c_j,r_j-s_n)\Big|\ \le\ \#\Big(\bigcup_j B(c_j,r_j)\Big)_{\Z}\ \le\ \Big|\bigcup_j B(c_j,r_j+s_n)\Big| .
\]
\end{lemma}

\begin{proof}
The cubes $Q_w\coloneqq w+[-\frac12,\frac12)^n$, $w\in\Z^n$, partition $\R^n$, have volume $1$, and satisfy $|x-w|\le s_n$ for $x\in Q_w$. If $w\in B(c_j,r_j)$ then $Q_w\subset B(c_j,r_j+s_n)$, which gives the upper bound. If $x\in B(c_j,r_j-s_n)\cap Q_w$, then $|w-c_j|<r_j$; hence the left-hand union is covered by the cubes $Q_w$ with $w$ in the middle set, which gives the lower bound.
\end{proof}

\begin{lemma}\label{lem:halo}
Let $\Et\subset\Z^n$ be finite and nonempty and $0<\alpha<1$.
\begin{enumerate}[label=\textup{(\roman*)},leftmargin=2.2em]
\item $m\in\Et_\alpha$ if and only if there is an open ball $B\ni m$ with $\#(B\cap\Et)>\alpha\,\#\Lat{B}$.
\item $\Et\subset\Et_\alpha$.
\item $\Et_\alpha$ is finite.
\end{enumerate}
\end{lemma}

\begin{proof}
(i) For $B\in\BHL$ we have $\Lat{m+B}=m+\Lat{B}$, and $B\mapsto m+B$ is a bijection from $\BHL$ onto the set of open balls containing $m$. (ii) If $m\in\Et$, the ball $B(m,1/2)$ contains no lattice point other than $m$. (iii) Let $m\in\Et_\alpha$ and let $B=B(c,r)$ be as in (i). Then $\#\Lat{B}<\#\Et/\alpha$, so Lemma~\ref{lem:count} gives $\kappa_n(r-s_n)^n<\#\Et/\alpha$, i.e.\ $r<r_0\coloneqq s_n+(\#\Et/(\alpha\kappa_n))^{1/n}$. Since $B\cap\Et\ne\emptyset$, $m$ lies within distance $2r_0$ of $\Et$.
\end{proof}

The next inequality follows from a theorem of Csik\'os~\cite[Theorem~4.2]{Cs}; an elementary proof is given in \cite[Lemma~3.1]{Mu}.

\begin{lemma}\label{lem:dilation}
Let $B_1,\dots,B_N$ be balls in $\R^n$ and $a\ge1$. Then $\big|\bigcup_jaB_j\big|\le a^n\big|\bigcup_jB_j\big|$.
\end{lemma}

The key geometric input is the following covering theorem.

\begin{theorem}[{Mukherjee~\cite[Theorem~A]{Mu}}]\label{thm:covering}
Let $n\ge2$ and $0<\delta<1/2$. Every finite family $\mathcal F$ of balls in $\R^n$ has a subfamily $\mathcal S$ such that
\[
\bigcup_{B\in\mathcal F}B\subset\bigcup_{B\in\mathcal S}(1+\delta)B\qquad\text{and}\qquad\sum_{B\in\mathcal S}\one_B\le A_n\,\delta^{-(n-1)/2},
\]
where $A_n$ depends only on $n$.
\end{theorem}

The exponent $(n-1)/2$ is optimal \cite{Mu}. Earlier, Chanillo and Muckenhoupt~\cite{CM} obtained the multiplicity bound $4^n\delta^{-n}$, and Han and Lu~\cite{HL} obtained $C_n\delta^{-n/2}\log(1/\delta)$.

Finally, we need the multidimensional Rokhlin lemma.

\begin{lemma}[{Katznelson--Weiss~\cite{KW}}]\label{lem:rokhlin}
Let $U_1,\dots,U_n$ be commuting invertible measure-preserving transformations of a standard probability space that generate a free action of $\Z^n$. For every $K\in\N$ and $\eta>0$ there is a measurable set $\Omega_0$ such that the sets $U^j\Omega_0$, $j\in\{0,\dots,K-1\}^n$, are pairwise disjoint and their union has measure greater than $1-\eta$.
\end{lemma}

We shall also use the following elementary lattice point count, which is a form of Riemann sum approximation.

\begin{lemma}\label{lem:jordan}
Let $A\subset\R^n$ be bounded with $|\partial A|=0$. Then $\varepsilon^n\#\Lat{\varepsilon^{-1}A}\to|A|$ as $\varepsilon\to0$.
\end{lemma}

\begin{proof}
Note that $\#\Lat{\varepsilon^{-1}A}=\#\{m\in\Z^n:\varepsilon m\in A\}$. Let $Q_m\coloneqq\varepsilon(m+[0,1)^n)$, and let $I\coloneqq\{m:Q_m\subset A\}$ and $O\coloneqq\{m:\overline{Q_m}\cap\overline A\ne\emptyset\}$. Since $\varepsilon m\in Q_m$, we have $I\subset\{m:\varepsilon m\in A\}\subset O$; moreover $\bigcup_{m\in I}Q_m\subset A\subset\bigcup_{m\in O}Q_m$. Hence both $\varepsilon^n\#\Lat{\varepsilon^{-1}A}$ and $|A|$ lie in the interval $[\varepsilon^n\#I,\varepsilon^n\#O]$. If $m\in O\setminus I$, then $\overline{Q_m}$ meets $\partial A$: otherwise the connected set $\overline{Q_m}$ would lie in one of the disjoint open sets $\operatorname{int}A$ and $\R^n\setminus\overline A$, and since it meets $\overline A$ we would have $\overline{Q_m}\subset\operatorname{int}A$, contradicting $m\notin I$. Thus the cells $Q_m$, $m\in O\setminus I$, are pairwise disjoint subsets of the closed $\sqrt n\varepsilon$-neighbourhood $(\partial A)_{\sqrt n\varepsilon}$ of the compact set $\partial A$, and
\[
\varepsilon^n(\#O-\#I)\le\big|(\partial A)_{\sqrt n\varepsilon}\big|\to|\partial A|=0\qquad(\varepsilon\to0). \qedhere
\]
\end{proof}

\section{The upper bound for the discrete constant}\label{sec:upper}

\begin{proposition}\label{prop:upper}
Let $n\ge2$ and $1/2\le\alpha<1$, and set $r_*\coloneqq(\kappa_n(1-\alpha))^{-1/n}-s_n$. If $r_*\ge4s_n$, then for every $0<\delta<1/2$,
\begin{equation}\label{eq:master}
\CtHL(\alpha)\le(1+\delta)^n\Big(1+\frac{3s_n}{r_*}\Big)^n\Big(1+A_n\delta^{-\frac{n-1}{2}}\,\frac{1-\alpha}{\alpha}\Big),
\end{equation}
where $A_n$ is the constant of Theorem~\ref{thm:covering}. Consequently there are $C_n>0$ and $\alpha_n\in(0,1)$ such that $\CtHL(\alpha)-1\le C_n(1-\alpha)^{1/n}$ for $\alpha_n<\alpha<1$.
\end{proposition}

\begin{proof}
Fix a finite nonempty set $\Et\subset\Z^n$. By Lemma~\ref{lem:halo}, $\Et\subset\Et_\alpha$ and $\Et_\alpha$ is finite.

\emph{Step 1: bad balls are large.} For each $m\in\Et_\alpha\setminus\Et$, Lemma~\ref{lem:halo}(i) provides an open ball $B_m=B(c_m,r_m)\ni m$ with $\#(B_m\cap\Et)>\alpha\,\#\Lat{B_m}$. Since $m\in\Lat{B_m}\setminus\Et$,
\begin{equation}\label{eq:holes}
1\le\#\big(\Lat{B_m}\setminus\Et\big)<(1-\alpha)\,\#\Lat{B_m}<\frac{1-\alpha}{\alpha}\,\#(B_m\cap\Et).
\end{equation}
In particular $\#\Lat{B_m}>(1-\alpha)^{-1}$, so Lemma~\ref{lem:count} gives $\kappa_n(r_m+s_n)^n>(1-\alpha)^{-1}$, that is, $r_m>r_*$.

\emph{Step 2: covering.} Apply Theorem~\ref{thm:covering} to the finite family $\{B_m\}_{m\in\Et_\alpha\setminus\Et}$, obtaining a subfamily $\mathcal S$, and let
\[
U\coloneqq\bigcup_{B\in\mathcal S}B,\qquad V\coloneqq\bigcup_{B\in\mathcal S}(1+\delta)B .
\]
Then $\Et_\alpha\subset\Et\cup\Lat{V}$, and since $U\subset V$,
\begin{equation}\label{eq:split}
\#\Et_\alpha\le\#\big(\Et\cup\Lat{U}\big)+\big(\#\Lat{V}-\#\Lat{U}\big).
\end{equation}

\emph{Step 3: holes in $U$.} By \eqref{eq:holes} and the multiplicity bound,
\[
\#\big(\Lat{U}\setminus\Et\big)\le\sum_{B\in\mathcal S}\#\big(\Lat{B}\setminus\Et\big)<\frac{1-\alpha}{\alpha}\sum_{w\in\Et}\sum_{B\in\mathcal S}\one_B(w)\le K\,\#\Et,\qquad K\coloneqq A_n\delta^{-\frac{n-1}{2}}\frac{1-\alpha}{\alpha}.
\]
Hence both $\#(\Et\cup\Lat{U})$ and $\#\Lat{U}$ are at most $(1+K)\#\Et$.

\emph{Step 4: dilation on the lattice.} Let $a\coloneqq(1+\delta)(1+3s_n/r_*)$. For every $B(c,r)\in\mathcal S$ we have $r>r_*\ge4s_n$, hence
\[
a(r-s_n)\ge(1+\delta)\Big(1+\frac{3s_n}{r}\Big)(r-s_n)=(1+\delta)\Big(r+2s_n-\frac{3s_n^2}{r}\Big)\ge(1+\delta)r+s_n ,
\]
so that $B(c,(1+\delta)r+s_n)\subset aB(c,r-s_n)$. By Lemma~\ref{lem:count}, Lemma~\ref{lem:dilation} and Lemma~\ref{lem:count} again,
\[
\#\Lat{V}\le\Big|\bigcup_{\mathcal S}B\big(c,(1+\delta)r+s_n\big)\Big|\le\Big|\bigcup_{\mathcal S}aB(c,r-s_n)\Big|\le a^n\Big|\bigcup_{\mathcal S}B(c,r-s_n)\Big|\le a^n\,\#\Lat{U}.
\]
Thus $\#\Lat{V}-\#\Lat{U}\le(a^n-1)(1+K)\#\Et$.

\emph{Step 5: conclusion.} By \eqref{eq:split} and Steps~3--4, $\#\Et_\alpha\le a^n(1+K)\#\Et$, which is \eqref{eq:master}.

\emph{Step 6: choice of parameters.} Let $t\coloneqq(1-\alpha)^{1/n}$, $\delta\coloneqq t$, and $X\coloneqq(\kappa_n(1-\alpha))^{-1/n}=\kappa_n^{-1/n}t^{-1}$, so that $r_*=X-s_n$. Choose $\alpha_n\in[1/2,1)$ so that for $\alpha>\alpha_n$ we have $t<1/2$ and $X\ge5s_n$. Then $r_*=X-s_n\ge4s_n$, so \eqref{eq:master} applies. Moreover, since $s_n\le X/5$, we have $r_*\ge\frac45X$, and therefore
\[
\frac{3s_n}{r_*}\le\frac{15s_n}{4X}=\lambda_nt,\qquad\lambda_n\coloneqq\tfrac{15}{4}s_n\kappa_n^{1/n}.
\]
Hence $a\le(1+t)(1+\lambda_nt)\le1+2(1+\lambda_n)t$ and, after enlarging $\alpha_n$ if necessary, $a^n\le1+C_nt$. On the other hand, since $1-\alpha=t^n$ and $\alpha\ge1/2$,
\[
K\le2A_n\,t^{\,n-\frac{n-1}{2}}=2A_n\,t^{\frac{n+1}{2}}\le2A_n\,t .
\]
Inserting these bounds into \eqref{eq:master} gives $\CtHL(\alpha)\le(1+C_nt)(1+2A_nt)\le1+C_n't$.
\end{proof}

\begin{remark}\label{rem:n1}
For $n=1$ we have $\BHL=\BS$, and the upper bound $\CtHL(\alpha)-1\lesssim1-\alpha$ follows from \cite[Section~4]{HP16} together with the one-dimensional case of \cite{HP14}. Alternatively, the classical fact that every finite family of intervals has a subfamily with the same union and overlap at most $2$ allows one to run Steps~2--5 with $V=U$, giving $\CtHL(\alpha)\le1+2(1-\alpha)/\alpha$.
\end{remark}

\section{The lower bound for the discrete constant}\label{sec:lower}

\begin{proposition}\label{prop:lower}
Let $n\ge1$ and $0<\alpha<1$ with $(\kappa_n(1-\alpha))^{-1/n}\ge\sqrt n$. Then
\[
\CtHL(\alpha)-1\ge\frac{\kappa_n^{1/n}}{2\cdot5^n}\,(1-\alpha)^{1/n}.
\]
\end{proposition}

\begin{proof}
Let $X\coloneqq(\kappa_n(1-\alpha))^{-1/n}$ and $\rho\coloneqq X+\sqrt n$, so $1\le\rho\le2X$. By Lemma~\ref{lem:count}, every ball $B$ of radius $\rho$ satisfies $\#\Lat{B}\ge\kappa_n(\rho-s_n)^n>\kappa_nX^n=(1-\alpha)^{-1}$. Let $R\coloneqq\lceil4\rho\rceil$ and $\Et\coloneqq\{0,1,\dots,R-1\}^n$, and fix $0<\varepsilon<1/(2\rho)$.

For $p'=(p_1,\dots,p_{n-1})\in\Z^{n-1}$ with $\rho\le p_i\le R-\rho$ for all $i$ (no condition if $n=1$), let
\[
p\coloneqq(p',-1),\qquad c\coloneqq p+(\rho-\varepsilon)e_n,\qquad B\coloneqq B(c,\rho).
\]
Then $p\in B$. We claim that $\Lat{B}\setminus\Et=\{p\}$. Let $w\in\Lat{B}$. Since $c_n=\rho-1-\varepsilon$, we have $-1-\varepsilon<w_n<2\rho-1-\varepsilon<R$; as $w_n$ is an integer and $\varepsilon<1$, this gives $-1\le w_n\le R-1$. If $w_n=-1$, then $|w'-p'|^2=|w-c|^2-(\rho-\varepsilon)^2<2\rho\varepsilon<1$, so $w=p$. If $w_n\ge0$, then $|w_i-p_i|<\rho$ for $i\le n-1$, whence $0\le p_i-\rho<w_i<p_i+\rho\le R$; thus $w\in\Et$. This proves the claim, and therefore
\[
\#(B\cap\Et)=\#\Lat{B}-1>\alpha\,\#\Lat{B},
\]
the inequality being equivalent to $\#\Lat{B}>(1-\alpha)^{-1}$, which holds by the choice of $\rho$ in the first paragraph (for every ball of radius $\rho$, wherever its center lies). Since $p\in B$, Lemma~\ref{lem:halo}(i) gives $p\in\Et_\alpha\setminus\Et$; the center of $B$ plays no role, because the maximal operator is uncentered.

Each $p_i$ ranges over the integers of a closed interval of length $R-2\rho\ge2\rho$, of which there are at least $\lfloor2\rho\rfloor\ge\rho$. Hence $\#(\Et_\alpha\setminus\Et)\ge\rho^{n-1}$, and since $R\le5\rho$,
\[
\CtHL(\alpha)-1\ge\frac{\rho^{n-1}}{R^n}\ge\frac{1}{5^n\rho}\ge\frac{1}{2\cdot5^nX}=\frac{\kappa_n^{1/n}}{2\cdot5^n}(1-\alpha)^{1/n}. \qedhere
\]
\end{proof}

\begin{proof}[Proof of Theorem~\ref{thm:discrete}]
Combine Proposition~\ref{prop:upper} (or Remark~\ref{rem:n1} if $n=1$) with Proposition~\ref{prop:lower}. \qedhere
\end{proof}

\section{Two-sided transference}\label{sec:transference}

\begin{lemma}\label{lem:freepart}
Let $U_1,\dots,U_n$ be commuting invertible measure-preserving transformations of a standard probability space $(\Omega,\Sigma,\mu)$, and let $F\coloneqq\{\omega:U^j\omega\ne\omega\text{ for all }j\ne0\}$. Then $F$ is measurable and invariant. If $\mu(F)>0$, the action admits arbitrarily large Rokhlin towers.
\end{lemma}

\begin{proof}
Since $(\Omega,\Sigma)$ is standard, there is a countable family $\{A_k\}_{k\ge1}\subset\Sigma$ separating the points of $\Omega$ (for instance, a countable base of a compatible Polish topology). For every $j\in\Z^n$,
\[
\{\omega:U^j\omega=\omega\}=\bigcap_{k\ge1}\big(\Omega\setminus(A_k\triangle U^{-j}A_k)\big).
\]
Indeed, if $U^j\omega=\omega$, then for every $k$ we have $\omega\in A_k$ if and only if $\omega\in U^{-j}A_k$; if $U^j\omega\ne\omega$, a set $A_k$ separating $\omega$ from $U^j\omega$ satisfies $\omega\in A_k\triangle U^{-j}A_k$. Hence each set $\{\omega:U^j\omega=\omega\}$ is measurable, and so is $F$, whose complement is the countable union of these sets over $j\ne0$. (This is where the standardness of the space is used.) Since $U^jU^k\omega=U^k\omega$ if and only if $U^j\omega=\omega$, we have $U^kF=F$ for all $k$. If $\mu(F)>0$, the action restricted to the standard probability space $(F,\Sigma|_F,\mu/\mu(F))$ is free, and Lemma~\ref{lem:rokhlin} (with $\eta=1/2$) produces, for each $K$, a set $\Omega_0\subset F$ of positive measure whose translates $U^j\Omega_0$, $j\in\{0,\dots,K-1\}^n$, are pairwise disjoint.
\end{proof}

\begin{proof}[Proof of Theorem~\ref{thm:transference}]
Part (i) is \cite[Theorem~3.1]{HP16}, and the final assertion follows from Lemma~\ref{lem:freepart}. Under the hypothesis of (ii), it remains to prove that $C^*_\cB(\alpha)\ge\Ct_\cB(\alpha)$.

Fix a finite nonempty set $\Et\subset\Z^n$ and a finite set $F\subset\{m:\Mt_\cB\one_{\Et}(m)>\alpha\}$. For each $m\in F$ choose $R_m\in\cB$ with
\[
\#\big((m+R_m)\cap\Et\big)>\alpha\,\#\big(R_m\cap\Z^n\big).
\]
Since the sets $R_m$ are bounded, there is $L\in\N$ such that $\Et$, $F$ and all the sets $(m+R_m)\cap\Z^n$ are contained in $[-L,L]^n$. Let $s\coloneqq(L,\dots,L)$, apply the hypothesis with $K=2L+1$ to obtain $\Omega_0$, and define
\[
E\coloneqq\bigsqcup_{k\in\Et}U^{k+s}\Omega_0,\qquad\text{so that}\qquad\mu(E)=\#\Et\,\mu(\Omega_0).
\]
Let $m\in F$ and $\omega=U^{m+s}\omega_0$ with $\omega_0\in\Omega_0$. For $j\in R_m\cap\Z^n$ we have $U^j\omega=U^{j+m+s}\omega_0$ and $j+m+s\in\{0,\dots,2L\}^n$. Since the levels $U^i\Omega_0$, $i\in\{0,\dots,2L\}^n$, are pairwise disjoint and $E$ is the union of the levels with indices in $\Et+s\subset\{0,\dots,2L\}^n$, we get $U^j\omega\in E$ if and only if $j+m\in\Et$. Consequently
\[
M^*_\cB\one_E(\omega)\ge\frac{1}{\#(R_m\cap\Z^n)}\sum_{j\in R_m\cap\Z^n}\one_E(U^j\omega)=\frac{\#((m+R_m)\cap\Et)}{\#(R_m\cap\Z^n)}>\alpha .
\]
Thus $M^*_\cB\one_E>\alpha$ on each of the pairwise disjoint sets $U^{m+s}\Omega_0$, $m\in F$, and
\[
\mu(\{M^*_\cB\one_E>\alpha\})\ge\#F\,\mu(\Omega_0)=\frac{\#F}{\#\Et}\,\mu(E).
\]
Taking the supremum over finite $F$ and then over $\Et$ gives $C^*_\cB(\alpha)\ge\Ct_\cB(\alpha)$.
\end{proof}

\begin{proof}[Proof of Corollary~\ref{cor:ergodic}]
By Theorem~\ref{thm:transference}(i) and Theorem~\ref{thm:discrete}, \mbox{$\CsHL(\alpha)-1\le\CtHL(\alpha)-1\le C_n(1-\alpha)^{1/n}$} for every action. If the action admits arbitrarily large Rokhlin towers, Theorem~\ref{thm:transference}(ii) and Theorem~\ref{thm:discrete} give $\CsHL(\alpha)-1=\CtHL(\alpha)-1\ge c_n(1-\alpha)^{1/n}$. Free actions on $\T$ exist for every $n$ (see Example~\ref{ex:examples}(a)), so no exponent larger than $1/n$ can hold uniformly.
\end{proof}

The following examples show that it is the size of the available towers, and not ergodicity, that determines whether the worst case is attained.

\begin{example}\label{ex:examples}
\begin{enumerate}[label=(\alph*),leftmargin=2em]
\item \emph{Free but not ergodic.} Let $\Omega=\T\times[0,1]$ with Lebesgue measure and $U_i(x,y)=(x+\theta_i,y)$, where $1,\theta_1,\dots,\theta_n$ are linearly independent over $\mathbb Q$. For $j\ne0$, $U^j$ translates the first coordinate by $\sum_ij_i\theta_i\notin\Z$, so the action is free; the set $\T\times[0,\frac12]$ is invariant, so it is not ergodic. By Theorem~\ref{thm:transference}, $C^*_\cB=\Ct_\cB$ for every basis $\cB$.

\item \emph{Ergodic with bounded periods.} Let $\Omega=(\Z/N\Z)^n$ with uniform measure and let $U_i$ be the cyclic shift by $e_i$. The action is transitive, hence ergodic. We claim that $\CsHL(\alpha)=1$ for $1-c\le\alpha<1$, where $c\coloneqq\min\{(3/5)^nN^{-n},(\kappa_n(5Ns_n)^n)^{-1}\}$.

First, every open ball $B=B(x,r)$ containing the origin satisfies $\#(B\cap N\Z^n)\ge c\,\#\Lat{B}$. Indeed, if $r\ge4Ns_n$, Lemma~\ref{lem:count} applied to the lattice $N\Z^n$ gives \\
$\#(B\cap N\Z^n)\ge\kappa_nN^{-n}(r-Ns_n)^n\ge\kappa_nN^{-n}(3r/4)^n$, while \mbox{$\#\Lat{B}\le\kappa_n(r+s_n)^n\le\kappa_n(5r/4)^n$}. If $r<4Ns_n$, then $0\in B\cap N\Z^n$ and $\#\Lat{B}\le\kappa_n(5Ns_n)^n$.

Now let $\emptyset\ne E\subset\Omega$ and $\omega\notin E$. Every $j\in B\cap N\Z^n$ satisfies $U^j\omega=\omega\notin E$, so $\MsHL\one_E(\omega)\le1-c$. Hence $\{\MsHL\one_E>\alpha\}=E$ for $\alpha\ge1-c$ (the inclusion \mbox{$E\subset\{\MsHL\one_E>\alpha\}$} follows from the ball $B(0,\frac12)$), and $\CsHL(\alpha)=1$, whereas $\CtHL(\alpha)>1$ for $\alpha$ close to $1$ by Theorem~\ref{thm:discrete}.

\item \emph{Nowhere free, with arbitrarily large towers.} Let $N_k\to\infty$ and $w_k>0$ with $\sum_kw_k=1$. Let $\Omega=\bigsqcup_k(\Z/N_k\Z)^n$, where the $k$-th component carries the uniform measure of total mass $w_k$, and let $U_i$ act on each component by the cyclic shift by $e_i$. Every point has nontrivial stabilizer, but for $K\le N_k$ the singleton $\Omega_0=\{0\}$ in the $k$-th component has positive measure and pairwise disjoint translates $U^j\Omega_0$, $j\in\{0,\dots,K-1\}^n$. Hence $C^*_\cB=\Ct_\cB$ for every basis $\cB$.
\end{enumerate}
\end{example}

\begin{remark}\label{rem:HS}
Theorem~\ref{thm:transference}(ii) is in the spirit of the transference principle of Hagelstein and Stokolos~\cite[Theorem~1(ii)]{HS}, who showed for non-periodic actions that weak type bounds for multiparameter ergodic maximal operators imply weak type bounds for the corresponding geometric operators, by painting lattice data onto the levels of a Katznelson--Weiss tower. Their argument transfers to $\R^n$ and loses dimensional constants, which is harmless for weak type bounds but not for Solyanik estimates. Transferring to $\Z^n$ instead, and working with indicator functions only, makes the transference exact.
\end{remark}

\section{Discrete versus Euclidean constants}\label{sec:comparison}

Both $\BHL$ and $\BS$ have the following properties: every element is a bounded open convex set; the basis is invariant under dilations $R\mapsto\lambda R$, $\lambda>0$; and if $R\in\cB$ and $z\in R$, then $R-z\in\cB$.

\begin{proposition}\label{prop:comparison}
Let $\cB\in\{\BHL,\BS\}$ and $0<\alpha<\beta<1$. Then $C_\cB(\beta)\le\Ct_\cB(\alpha)$.
\end{proposition}

\begin{proof}
Fix a measurable set $E$ with $0<|E|<\infty$ and let $H\coloneqq\{M_\cB\one_E>\beta\}$. Call an open set $P$ \emph{good} if $P=x+R$ for some $x\in\R^n$ and $R\in\cB$, and $|P\cap E|>\beta|P|$. By the last property of $\cB$, if $P$ is good then $P\subset H$; hence $H$ is the union of all good sets and is open. Let $K\subset H$ be compact, and cover $K$ by finitely many good sets $P_1,\dots,P_N$.

Let $\eta>0$. By the regularity of Lebesgue measure, there exist $k\in\N$ and a finite union $G$ of half-open dyadic cubes of side length $2^{-k}$ such that $|E\triangle G|<\min\{\eta,\frac{\beta-\alpha}{2}\min_i|P_i|\}$; here $\min_i|P_i|>0$, since the $P_i$ are finitely many nonempty open sets. Since $P_i\cap E\subset(P_i\cap G)\cup(E\setminus G)$ and $P_i$ is good, for every $i$
\[
|P_i\cap G|\ge|P_i\cap E|-|E\triangle G|>\beta|P_i|-\tfrac{\beta-\alpha}{2}|P_i|=\tfrac{\alpha+\beta}{2}|P_i| .
\]
For $\varepsilon=2^{-\ell}$ with $\ell\ge k$, let $\Et_\varepsilon\coloneqq\{m\in\Z^n:\varepsilon m\in G\}$. Every dyadic cell $\varepsilon(m+[0,1)^n)$ is either contained in $G$ or disjoint from it, so $G=\bigcup_{m\in\Et_\varepsilon}\varepsilon(m+[0,1)^n)$ and $\varepsilon^n\#\Et_\varepsilon=|G|$.

The sets $P_i$, $P_i\cap G$ and $\bigcup_iP_i$ are bounded and have boundaries of measure zero: $\partial P_i$ is a sphere or the boundary of a rectangle, $\partial G$ is contained in finitely many hyperplanes, and $\partial(A\cap A')\cup\partial(A\cup A')\subset\partial A\cup\partial A'$. Since $(\varepsilon^{-1}P_i)\cap\Et_\varepsilon=\{m\in\Z^n:\varepsilon m\in P_i\cap G\}=\Lat{\varepsilon^{-1}(P_i\cap G)}$, Lemma~\ref{lem:jordan} gives, as $\varepsilon\to0$,
\[
\varepsilon^n\#\Lat{\varepsilon^{-1}P_i}\to|P_i|,\qquad\varepsilon^n\#\big(\varepsilon^{-1}P_i\cap\Et_\varepsilon\big)\to|P_i\cap G|,\qquad\varepsilon^n\#\Lat{\varepsilon^{-1}\textstyle\bigcup_iP_i}\to\big|\textstyle\bigcup_iP_i\big| .
\]
Hence for all sufficiently small $\varepsilon$ we have $\#(\varepsilon^{-1}P_i\cap\Et_\varepsilon)>\alpha\,\#\Lat{\varepsilon^{-1}P_i}$ for every $i$. By dilation invariance, $\varepsilon^{-1}P_i$ is a translate of an element of $\cB$, and by the last property of $\cB$, $\varepsilon^{-1}P_i-m\in\cB$ for every $m\in\Lat{\varepsilon^{-1}P_i}$. Therefore every lattice point of $\varepsilon^{-1}\bigcup_iP_i$ lies in $\{\Mt_\cB\one_{\Et_\varepsilon}>\alpha\}$, and
\[
\Ct_\cB(\alpha)\ge\liminf_{\varepsilon\to0}\frac{\#\Lat{\varepsilon^{-1}\bigcup_iP_i}}{\#\Et_\varepsilon}\ge\frac{|K|}{|G|}\ge\frac{|K|}{|E|+\eta}.
\]
Letting $\eta\to0$ and taking the supremum over compact $K\subset H$ and over $E$ gives $C_\cB(\beta)\le\Ct_\cB(\alpha)$.
\end{proof}

\begin{proof}[Proof of Theorem~\ref{thm:comparison}]
By \cite{HP15} (see also \cite[Problem~7.4]{HP16}), $\CHL$ and $\CS$ are continuous on $(0,1)$. Letting $\beta\downarrow\alpha$ in Proposition~\ref{prop:comparison} gives $C_\cB(\alpha)\le\Ct_\cB(\alpha)$.
\end{proof}

\begin{proof}[Proof of Corollary~\ref{cor:strong}]
The reverse inequality $\CtS(\alpha)\le\CS(\alpha)$, $0<\alpha<1$, is proved in \cite[Section~4]{HP16}; we recall the argument. Given a finite set $\Et\subset\Z^n$, let $E\coloneqq\Et+[0,1)^n$, so that $|E|=\#\Et$. For $R\in\BS$, the unit cells $j+[0,1)^n$, $j\in R\cap\Z^n$, tile an axis-parallel rectangle containing $[0,1)^n$. Hence every lattice average of $\one_{\Et}$ over a translate of $R$ is a Euclidean average of $\one_E$ over an axis-parallel rectangle, so that $M_{\BS}\one_E\ge\Mt_{\BS}\one_{\Et}(m)$ on the open cell $m+(0,1)^n$, and
\[
\#\{\Mt_{\BS}\one_{\Et}>\alpha\}\le\big|\{M_{\BS}\one_E>\alpha\}\big|\le\CS(\alpha)\,\#\Et .
\]
Combined with Theorem~\ref{thm:comparison}, this gives $\CtS=\CS$. The remaining assertions follow from Theorem~\ref{thm:transference}(ii) and the last statement of Theorem~\ref{thm:transference}.
\end{proof}

\begin{remark}\label{rem:strongHolder}
Corollary~\ref{cor:strong} shows that for free actions the ergodic constant $\CsS$ coincides with the Euclidean halo function $\CS$. In particular $\CsS$ is locally H\"older continuous on $(0,1)$ with the exponent obtained in \cite{HP15}, and its Solyanik behaviour is exactly that of $\CS$. Whether $\CtHL$, and hence $\CsHL$ for free actions, is continuous on $(0,1)$ remains open.
\end{remark}

The inequality $\CHL\le\CtHL$ of Theorem~\ref{thm:comparison} is strict for $\alpha$ close to $1$ when $n\ge2$. The following example makes this explicit.

\begin{example}\label{ex:strict}
Let $n=2$. We show that $\CHL(\alpha)<\CtHL(\alpha)$ for every $\alpha$ with $0<1-\alpha\le10^{-7}$. Write $s\coloneqq1-\alpha$.

\emph{Discrete lower bound.} Let $X\coloneqq(\pi s)^{-1/2}$ (so $X\ge1784$), $\rho\coloneqq X+1$, $H\coloneqq\lceil2\rho-1\rceil$, and for $W\in\N$ let $\Et_W\coloneqq\{0,\dots,W-1\}\times\{0,\dots,H-1\}$. By Lemma~\ref{lem:count}, every disc of radius $\rho$ contains at least $\pi(\rho-\frac{\sqrt2}{2})^2>\pi X^2=s^{-1}$ lattice points. Fix $0<\varepsilon<1/(2\rho)$. For each integer $p_1$ with $\rho-1\le p_1\le W-\rho$, let $p\coloneqq(p_1,-1)$ and $B\coloneqq B((p_1,\rho-1-\varepsilon),\rho)$. Then $p\in B$, and exactly as in the proof of Proposition~\ref{prop:lower} one checks that $\Lat{B}\setminus\Et_W=\{p\}$: every $w\in\Lat B$ satisfies $-1-\varepsilon<w_2<2\rho-1-\varepsilon<H$ and $p_1-\rho<w_1<p_1+\rho$, so that $0\le w_1\le W-1$ and $-1\le w_2\le H-1$, and $w_2=-1$ forces $|w_1-p_1|^2<2\rho\varepsilon<1$, i.e.\ $w=p$. Hence $\#(B\cap\Et_W)=\#\Lat B-1>\alpha\,\#\Lat B$, and $p\in(\Et_W)_\alpha$ by Lemma~\ref{lem:halo}(i). By the symmetry $(x,y)\mapsto(x,H-1-y)$ of $\Et_W$, the same holds for $(p_1,H)$. There are at least $W-2\rho$ admissible values of $p_1$, so
\[
\CtHL(\alpha)\ \ge\ 1+\frac{2(W-2\rho)}{WH}\ \xrightarrow[W\to\infty]{}\ 1+\frac2H\ \ge\ 1+\frac{1}{X+1}\ \ge\ 1+1.77\,s^{1/2},
\]
where we used $H\le2\rho=2X+2$ and $X\ge1784$.

\emph{Euclidean upper bound.} Mukherjee's estimate \cite[(3.3)]{Mu} (see also \eqref{eq:euclid} below) gives, for $0<\delta<1/2$,
\[
\CHL(\alpha)\le(1+\delta)^2\Big(1+A_2\,\delta^{-1/2}\,\frac{s}{\alpha}\Big),
\]
where $A_2$ is the constant of Theorem~\ref{thm:covering} for $n=2$. In the proof of \cite[Theorem~A]{Mu}, the selected discs containing a given point give rise to points on $S^1$ that are pairwise at distance at least $\sqrt{\delta/10}$. Consecutive such points are separated by an angle of at least $\sqrt{\delta/10}$, so there are at most $2\pi/\sqrt{\delta/10}$ of them, and one may take $A_2=2\pi\sqrt{10}<19.87$. Choosing $\delta=s^{2/3}\le2.2\cdot10^{-5}$ and using $s/\alpha\le1.0000002\,s$, we obtain
\[
\CHL(\alpha)\le\big(1+2.0001\,s^{2/3}\big)\big(1+19.871\,s^{2/3}\big)\le1+22\,s^{2/3}.
\]

\emph{Comparison.} Since $s^{1/6}\le10^{-7/6}<0.07<1.77/22$, we have $22\,s^{2/3}<1.77\,s^{1/2}$, and therefore
\[
\CHL(\alpha)\ \le\ 1+22\,s^{2/3}\ <\ 1+1.77\,s^{1/2}\ \le\ \CtHL(\alpha).
\]
For instance, for $\alpha=1-10^{-7}$ the finite set $\Et_W$ with $W=10^5$ (so that $H=3570$) already satisfies $\#(\Et_W)_\alpha/\#\Et_W>1.00054$, while optimizing in $\delta$ (for example $\delta=6\cdot10^{-5}$) gives $\CHL(\alpha)\le1.000377$.

The size of the threshold is dictated by the exponents rather than by the particular example: the two sides behave like $s^{1/2}$ and $s^{2/3}$, whose exponents differ by $1/6$, so the strict inequality can only be certified once $s^{1/6}$ is smaller than the ratio of the constants. The mechanism is the lattice first-layer effect: for a long strip $(0,W)\times(0,H)$ in $\R^2$ the Euclidean halo gains only a relative amount of order $s^{2/3}$, whereas on the lattice an almost tangent disc captures a whole extra row on each side, a relative gain of $2/H\asymp s^{1/2}$. We note that the constant $A_2$ is extracted from the proof in the preprint \cite{Mu}.
\end{example}

\section{Further remarks and questions}\label{sec:remarks}

\begin{remark}[Covering input]\label{rem:covering}
The proof of Proposition~\ref{prop:upper} only uses a covering theorem with multiplicity at most $A_n\delta^{-\beta}(\log(e/\delta))^{\gamma}$ for some $\beta<n-1$ and $\gamma\ge0$: with $\delta=t=(1-\alpha)^{1/n}$ one then has $K\lesssim t^{\,n-\beta}(\log(e/t))^{\gamma}\lesssim t$. The Han--Lu bound~\cite{HL} has $\beta=n/2$ and $\gamma=1$, which suffices for $n\ge3$. For $n=2$ it yields $\CtHL(\alpha)-1\lesssim(1-\alpha)^{1/2}\log\frac{1}{1-\alpha}$, and Theorem~\ref{thm:covering} removes the logarithm. The Chanillo--Muckenhoupt bound ($\beta=n$) is not sufficient.
\end{remark}

\begin{remark}[The gap between $1/n$ and $2/(n+1)$]
Mukherjee's proof of the upper bound in \eqref{eq:Mu} yields
\begin{equation}\label{eq:euclid}
\CHL(\alpha)\le(1+\delta)^n\Big(1+A_n\delta^{-\frac{n-1}{2}}\,\frac{1-\alpha}{\alpha}\Big),\qquad0<\delta<\tfrac12,
\end{equation}
and the choice $\delta=(1-\alpha)^{2/(n+1)}$ gives the exponent $2/(n+1)$. Inequality \eqref{eq:master} is precisely \eqref{eq:euclid} multiplied by the lattice factor $(1+3s_n/r_*)^n$. With $\delta=(1-\alpha)^{2/(n+1)}$, the $\delta$-dependent factors in \eqref{eq:master} are $1+O_n((1-\alpha)^{2/(n+1)})$, which is $o((1-\alpha)^{1/n})$ for $n\ge2$; the lattice factor contributes $3s_n/r_*\asymp_n(1-\alpha)^{1/n}$ and cannot be removed, as Proposition~\ref{prop:lower} shows. In the Euclidean lower bound the halo is produced by caps of height $h$, whose volume is of order $h^{(n+1)/2}$; on the lattice, a ball almost tangent to a lattice hyperplane captures a single lattice point of the next layer at no cost in radius, which is a stronger effect.
\end{remark}

\begin{remark}[The case $n=1$]
For $n=1$ all three constants $\CHL,\CtHL,\CsHL$ (the last for free actions) satisfy $C(\alpha)-1\asymp1-\alpha$, so all the exponents discussed above equal $1$. Indeed, $\BHL=\BS$ when $n=1$, so the upper bounds follow from \cite{HP14} and Remark~\ref{rem:n1}, while for $E=(0,1)$ and $0<t<1/\alpha-1$ the averages of $\one_E$ over $(-t,1)$ and $(0,1+t)$ exceed $\alpha$, which gives $\CHL(\alpha)\ge1+2(1/\alpha-1)$; the lower bounds for $\CtHL$ and $\CsHL$ follow from Proposition~\ref{prop:lower} and Theorem~\ref{thm:transference}.
\end{remark}

The contrast between balls (Example~\ref{ex:strict}) and rectangles (Corollary~\ref{cor:strong}) suggests the following question.

\begin{question}
Let $\mathcal K$ be a convex body containing the origin in its interior, and let $\cB_{\mathcal K}$ be the basis of all translates and dilates of $\mathcal K$ that contain the origin. How does the Solyanik exponent of $\Ct_{\cB_{\mathcal K}}$ compare with that of $C_{\cB_{\mathcal K}}$, and how does the answer depend on the curvature of $\partial\mathcal K$ and on the position of its faces relative to the lattice?
\end{question}

Example~\ref{ex:examples}(b) shows that $\CsHL$ can be much smaller than $\CtHL$ when the Rokhlin towers have bounded size. This leads to a quantitative form of \cite[Problem~7.4]{HP16}.

\begin{question}
For which actions that do not admit arbitrarily large Rokhlin towers is the Solyanik exponent of $\CsHL$ strictly larger than $1/n$? For instance, what is it for $n=2$, $U_1$ an irrational rotation of $\T$ and $U_2=\mathrm{Id}$, where each orbit is a copy of $\Z^2$ collapsed onto $\Z$?
\end{question}

\section*{Acknowledgements}
The author thanks Paul Hagelstein for suggesting this problem and for helpful discussions.

\end{document}